\documentclass[aos,MSNbibl,seceqn,citesort,dvips]{arximspdf}
\usepackage{graphicx}

\doi{10.1214/12-AOS969} %kopijuoti is PTS
\volume{40}
\issue{2}
\pubyear{2012}
\firstpage{666}
\lastpage{693}

\makeatletter

\newcommand{\dR}{\mathbb{R}}
\newcommand{\dE}{\mathbb{E}}
\newcommand{\dP}{\mathbb{P}}
\newcommand{\cG}{\mathcal{G}}
\newcommand{\cM}{\mathcal{M}}
\newcommand{\cN}{\mathcal{N}}
\newcommand{\cP}{\mathcal{P}}
\newcommand{\cQ}{\mathcal{Q}}
\newcommand{\cF}{\mathcal{F}}
\newcommand{\veps}{\varepsilon}
\newcommand{\rI}{\mathrm{I}}
\newcommand{\wh}{\widehat}

\newcommand{\liml}{\stackrel{\mathcal L}{\longrightarrow}}

\newtheorem{theorem}{Theorem}[section]
\newproclaim{rem}{Remark}[section]

\makeatother

\begin{document}
\begin{frontmatter}

\title{A Robbins--Monro procedure for estimation in semiparametric regression models}
\runtitle{A Robbins--Monro procedure}

\begin{aug}
\author[A]{\fnms{Bernard} \snm{Bercu}\corref{}\ead[label=e1]{Bernard.Bercu@math.u-bordeaux1.fr}}
\and
\author[A]{\fnms{Philippe} \snm{Fraysse}\ead[label=e2]{Philippe.Fraysse@math.u-bordeaux1.fr}}
\runauthor{B. Bercu and P. Fraysse}
\affiliation{Universit\'e de Bordeaux}
\address[A]{Institut de Math\'ematiques de Bordeaux\\
UMR CNRS\\
Universit\'e de Bordeaux\\
and \\
INRIA Bordeaux, Team ALEA\\
F-33400 Talence\\
France\\
\printead{e1}\\
\hphantom{E-mail: }\printead*{e2}} %adresu isvedimo komanda gale!
\end{aug}

% HISTORY:
\received{\smonth{2} \syear{2011}}
\revised{\smonth{1} \syear{2012}}

% ABSTRACT
%
\begin{abstract}
This paper is devoted to the parametric estimation of a shift together
with the nonparametric estimation of a regression function in a
semiparametric regression model. We implement a very efficient and easy
to handle Robbins--Monro procedure. On the one hand, we propose a
stochastic algorithm similar to that of Robbins--Monro in order to
estimate the shift parameter. A preliminary evaluation of the
regression function is not necessary to estimate the shift parameter.
On the other hand, we make use of a recursive Nadaraya--Watson
estimator for the estimation of the regression function. This kernel
estimator takes into account the previous estimation of the shift
parameter. We establish the almost sure convergence for both
Robbins--Monro and Nadaraya--Watson estimators. The asymptotic
normality of our estimates is also provided. Finally, we illustrate our
semiparametric estimation procedure on simulated and real data.
\end{abstract}

% KEYWORDS
%

\end{frontmatter}

%%%%%%%%%%%%%%%%%%%%%%%%%%%%%%%%%%%%%%%%%%%%%%%%%%%%%%%%%%%%%%%%%%%%%%%%%%%%%%%%%%%%

%s1 #&#
\section{Introduction}\label{sec1}

%%%%%%%%%%%%%%%%%%%%%%%%%%%%%%%%%%%%%%%%%%%%%%%%%%%%%%%%%%%%%%%%%%%%%%%%%%%%%%%%%%%%

A wide range of real-life phenomena occur periodically. One can think
about meteorology with daily or annual
cycles of temperature~\cite{MR2744122}, astronomy with the famous
11-year cycles of solar geomagnetic activity
\cite{Lassen}, medicine with human circadian rhythms~\cite{Wang} or
ECG signals~\cite{Trigano},
econometry~\cite{MR1623559,MR833929},
communication~\cite{MR2369028}, etc. Statistical analysis of periodic
data is of great interest in order
to design suitable models for those cyclic phenomena. An important
literature is available on the
so-called periodic shape-invariant model introduced by Lawton,
Sylvestre and Maggio~\cite{Lawton}.
Theoretical advances on shape-invariant models together with
statistical applications may be found in
\cite{MR1041386,MR2744122,MR1332581,MR924858,MR2662362,Wang}.
A periodic shape-invariant model is a semiparametric regression model with
an unknown periodic shape function. It is given, for all $n\geq0$, by
%%
%
%e1.1 #&#
%
\begin{equation}
\label{Sim}
Y_{n}=h(X_{n})+\veps_n,
\end{equation}
where the inputs $(X_n)$ are known observation times, the output
$(Y_n)$ are the observations,
and $(\veps_n)$ are unknown random errors. The function $h$ is
periodic and takes the form
\[
h(x)=m + \sum_{k=1}^p a_k f(x- \theta_k),
\]
where~$f$ represents the unknown characteristic shape function, $m$ is
the overall mean, while
$\theta=(\theta_1,\ldots,\theta_p)$ and $a=(a_1,\ldots,a_p)$ are
unknown shift and scale parameters.

In this paper, we shall focus our attention on the particular case
$p=1$, $m=0$ and $a=1$ by studying
the semiparametric regression model given, for all $n\geq0$, by
%%
%
%e1.2 #&#
%
\begin{equation}
\label{Sempar}
Y_{n}=f(X_{n}-\theta)+\veps_n,
\end{equation}
where $(X_n)$ and $(\veps_n)$ are two independent sequences of
independent and identically distributed random variables.
We are dealing with random observation times in contrast with the previous
literature where $(X_n)$ are assumed to be known and equidistributed
over a given interval.
We are interested in the parametric estimation of the shift parameter
$\theta$ together with the nonparametric estimation
of the shape function~$f$. However, one has to keep in mind that our
main interest lies in the estimation of the parameter $\theta$.
We are also motivated by a statistical application on the detection of
atrial fibrillation using ECG analysis
\cite{Clifford,Trigano}.

First of all, we implement a Robbins--Monro procedure in order to
estimate the unknown parameter $\theta$
without any preliminary evaluation of the regression function~$f$. Our
approach is very easy to handle and it performs
very well. Moreover, our approach is totally different from the one
recently proposed by Dalalyan, Golubev and
Tsybakov~\cite{MR2275239} in the Gaussian white-noise case. First, a
penalized maximum likelihood estimator
of $\theta$ is proposed in~\cite{MR2275239} with an appropriately
chosen penalty based on a Fourier series
approximation of the function~$f$. Second, the asymptotic behavior of
the mean squared risk of this estimator is investigated.
One can observe that our estimator is much easier to calculate. In
addition, we
do not require any assumption on the derivatives of the function~$f$.
In the situation where the
parameter $\theta$ is random, Castillo and Loubes~\cite{MR2508947}
propose a plug-in version of the Parzen--Rosenblatt
\cite{MR0143282,MR0343355} density estimator of $\theta$.
The construction of this estimate also relies on the
penalized maximum likelihood estimator of~$\theta$ given in
\cite{MR2275239}.
Furthermore, in the case where one observes several Gaussian functions
differing from
each other by a translation parameter, Gamboa, Loubes and Maza
\cite{MR2369028} propose
to transform the starting model by using a~discrete Fourier transform.
Hence, from the resulting model, they estimate the shift parameters by
minimizing a quadratic functional.
This approach is very interesting by the few assumptions made on the
regression function. In a more general framework,
Vimond~\cite{MR2662362} makes use of a truncated Fourier approximation
of~$f$ in order to evaluate the profile
log-likelihood score function associated with the shift and scale parameters.
This two-step strategy requires, as in~\cite{MR2369028}, the
estimation of the Fourier coefficients of~$f$.
However, it performs pretty well as it leads to consistent and
asymptotically efficient estimators
of the shift and scale parameters. Our alternative approach to estimate
$\theta$ is associated to a stochastic recursive
algorithm similar to that of Robbins--Monro
\cite{MR0042668,MR0343355}.

Assume that one can find a function $\phi$, free of the parameter
$\theta$, such that \mbox{$\phi(\theta)=0$}. Then, it is possible to estimate
$\theta$ by the Robbins--Monro algorithm
%
%e1.3 #&#
%
\begin{equation}
\label{RMalgo}
\widehat{\theta}_{n+1}=\widehat{\theta}_{n}+\gamma_{n}T_{n+1},
\end{equation}
where $(\gamma_n)$ is a positive sequence of real numbers decreasing
toward zero and~$(T_n)$ is a sequence of random variables such that
$\dE[T_{n+1}|\cF_n]=\phi(\widehat{\theta}_{n})$
where~$\mathcal{F}_{n}$ stands for the $\sigma$-algebra of the events
occurring up to time $n$.
Under standard conditions on the function $\phi$ and on the sequence
$(\gamma_n)$, it is well known~\cite{MR1485774,MR1993642} that
$\widehat{\theta}_{n}$ tends to $\theta$ almost surely. The
asymptotic normality of $\widehat{\theta}_{n}$ together
with the quadratic strong law may also be found in
\cite{MR624435,MR2351104} and~\cite{MR1654569}. A randomly truncated
version of the Robbins--Monro algorithm is also given in
\cite{MR931029,MR2542461}.

Our second goal is the estimation of the unknown regression function~$f$.
A wide range of literature is available on nonparametric estimation of
a regression function. We refer the reader to
\cite{MR1843146,MR2013911} for two excellent books on density
and regression function estimation.
Here, we focus our attention on the Nadaraya--Watson estimator of~$f$.
The almost sure convergence of the Nadaraya--Watson estimator
\cite{MR0166874,MR0185765} without the shift $\theta$ was
established by Noda
\cite{MR0426278}; see also H\"ardle et al.
\cite{MR740916,MR964932} for the law of iterated logarithm and the
uniform strong law.
A nice extension of the previous results may be found in
\cite{MR924860}. The asymptotic normality of the Nadaraya--Watson estimator
was proved by
Schuster~\cite{MR0301845}. Moreover, Choi, Hall and Rousson
\cite{MR1805786} propose three data-sharpening versions of the
Nadaraya--Watson estimator
in order to reduce the asymptotic variance in the central limit
theorem. Furthermore, in the situation where the regression function is
monotone,
Hall and Huang~\cite{MR1865334} provide a method for monotonizing the
Nadaraya--Watson estimator. For $n$ large enough, their alternative estimator
coincides with the standard Nadaraya--Watson estimator on a compact
interval where the regression function~$f$ is monotone.
In our situation, we propose to make use of a recursive
Nadaraya--Watson estimator~\cite{MR1485774} of~$f$ which
takes into account the previous estimation of the shift parameter
$\theta$.
It is given, for all $x\in\dR$, by
%
%e1.4 #&#
%
\begin{equation}
\label{RNW}
\widehat{f}_{n}(x)=\frac{\sum_{k=1}^{n} W_{k}(x)Y_{k}}{\sum
_{k=1}^{n} W_{k}(x)}
\end{equation}
with
\[
W_{n}(x)=\frac{1}{h_{n}}K\biggl(\frac{X_{n}-\widehat{\theta
}_{n-1}-x}{h_{n}}\biggr),
\]
where the kernel $K$ is a chosen probability density function and
the bandwidth $(h_n)$ is a sequence of positive real numbers
decreasing to zero. The main difficulty arising here is that we have to
deal with the additional term~$\widehat{\theta}_{n}$ inside the kernel~$K$. Consequently, we are
led to analyze a double stochastic algorithm
with, at the same time, the study of the asymptotic behavior of the
Robbins--Monro estimator
$\widehat{\theta}_{n}$ of~$\theta$, and the Nadaraya--Watson
estimator~$\widehat{f}_{n}$ of~$f$.

The paper is organized as follows. Section~\ref{sec2} is devoted to the
parametric estimation
of $\theta$. We establish the almost sure convergence of $\widehat
{\theta}_{n}$
as well as a~law of iterated logarithm and the asymptotic normality.
Section~\ref{sec3}
deals with the nonparametric estimation of~$f$. Under standard
regularity assumptions on the kernel~$K$,
we prove the almost sure pointwise convergence of~$\widehat{f}_{n}$ to~$f$. In addition,
we also establish the asymptotic normality of $\widehat{f}_{n}$.
Section~\ref{sec4} contains some numerical experiments on simulated and real ECG
data, illustrating
the performances of our semiparametric estimation procedure.
The proofs of the parametric results are given in Section~\ref{sec5}, while
those concerning
the nonparametric results are postponed to Section~\ref{sec6}.
%%%%%%%%%%%%%%%%%%%%%%%%%%%%%%%%%%%%%%%%%%%%%%%%%%%%%%%%%%%%%%%%%%%%%%%%%%%%%%%%%%%%

%s2 #&#
\section{Estimation of the shift}\label{sec2}

%%%%%%%%%%%%%%%%%%%%%%%%%%%%%%%%%%%%%%%%%%%%%%%%%%%%%%%%%%%%%%%%%%%%%%%%%%%%%%%%%%%%

First of all, we focus our attention on the estimation of the shift
parameter $\theta$ in the semiparametric regression model given
by~(\ref{Sempar}). We assume that $(\veps_n)$ is a sequence of
independent and identically distributed random variables
with zero mean and unknown positive variance $\sigma^2$. Moreover, it
is necessary to make several hypotheses similar to those of
\cite{MR2275239}.
{\renewcommand{\theequation}{$\mathcal{H}_{\arabic{equation}}$}
\begin{eqnarray}\label{equH1}\qquad
\\[-17pt]
&&\begin{tabular}{p{318pt}}
The observation times $(X_n)$ are
independent and identically distributed
with probability density function $g$, positive on its
support $[-1/2,1/2]$. In
addition, $g$ is continuous, twice differentiable with
bounded derivatives.
\end{tabular}\nonumber
\\
\label{equH2}
&&\begin{tabular}{p{300pt}}
The shape function~$f$ is symmetric,
bounded, periodic with period 1.
\end{tabular}
\end{eqnarray}}

\vspace*{-\baselineskip}

Let $X$ be a random variable sharing the same distribution as $(X_n)$.
In all the sequel, the auxiliary function $\phi$ defined, for all
$t\in\dR$, by
\setcounter{equation}{0}
%e2.1 #&#
%
\begin{equation}
\label{defphi}
\phi(t)=\dE\biggl[\frac{\sin(2\pi(X-t)) }{g(X)}f(X-\theta)\biggr]
\end{equation}
will play a prominent role. More precisely, it follows from the
periodicity of~$f$ that
\begin{eqnarray*}
\phi(t)&=&\int_{-1/2}^{1/2}\sin\bigl(2\pi(x-t)\bigr)f(x-\theta) \,dx\\
&=&\int
_{-1/2}^{1/2}\sin\bigl(2\pi(y+\theta-t)\bigr)f(y) \,dy,\\
&=&\sin\bigl(2\pi(\theta-t)\bigr)\int_{-1/2}^{1/2}\cos(2\pi y)f(y) \,dy\\
&&{}+\cos
\bigl(2\pi(\theta-t)\bigr)\int_{-1/2}^{1/2}\sin(2\pi y)f(y) \,dy.
\end{eqnarray*}
Consequently, the symmetry of~$f$ leads to
%
%e2.2 #&#
%
\begin{equation}
\label{DEFPHI}
\phi(t)=f_{1} \sin\bigl(2\pi(\theta-t)\bigr),
\end{equation}
where $f_1$ is the first Fourier coefficient of~$f$
\[
f_1=\int_{-1/2}^{1/2}\cos(2\pi x)f(x) \,dx.
\]
Throughout the paper, we assume that $f_1 \neq0$.
Obviously, $\phi$ is a continuous and bounded function such that $\phi
(\theta)=0$.
In addition, one can easily verify that for all $t\in\dR$
such that $|t-\theta|<1/2$, the product $(t-\theta)\phi(t)$
has a~constant sign. It is negative if $f_1>0$, while it is positive if $f_1<0$.
Therefore, we are in position to implement our Robbins--Monro procedure
\cite{MR0042668,MR0343355}.
Let $K=[-1/4,1/4]$ and denote by $\pi_{K}$ the projection on the
compact set $K$ defined, for all $x\in\dR$, by
\[
\pi_{K}(x) = \cases{
x, &\quad if $|x|\leq1/4$,\cr
1/4, &\quad if $x\geq1/4$,\cr
-1/4, &\quad if $x\leq-1/4$.}
\]
Let $(\gamma_n)$ be a decreasing sequence of positive real numbers satisfying
%
%e2.3 #&#
%
\begin{equation}
\label{hypgamma}
\sum_{n=1}^\infty\gamma_{n}=+\infty
\quad\mbox{and}\quad
\sum_{n=1}^\infty\gamma_{n}^2<+\infty.
\end{equation}
For the sake of clarity, we shall make use of $\gamma_n=1/n$.
We estimate the shift parameter $\theta$ via the projected
Robbins--Monro algorithm
%
%e2.4 #&#
%
\begin{equation}
\label{RMA}
\wh{\theta}_{n+1}=\pi_{K}\bigl(\wh{\theta}_{n}+
\operatorname{sign}(f_1)\gamma_{n+1}T_{n+1}\bigr),
\end{equation}
where the initial value $\wh{\theta}_{0} \in K$ and the random
variable $T_{n+1}$ is defined by
%
%e2.5 #&#
%
\begin{equation}
\label{DefT}
T_{n+1}=\frac{\sin(2\pi(X_{n+1}-\wh{\theta}_{n}))}{g(X_{n+1})}Y_{n+1}.
\end{equation}

Our first result concerns the almost sure convergence of the estimator
$\wh{\theta}_{n}$.
%
%th2.1 #&#
%
\begin{theorem}
\label{thmascvg}
Assume that~(\ref{equH1}) and~(\ref{equH2}) hold and that
$|\theta|<1/4$. Then, $\wh{\theta}_{n}$ converges almost surely to
$\theta$.
In addition, the number of times that
the random variable $\wh{\theta}_{n}+\operatorname{sign}(f_1)\gamma
_{n+1}T_{n+1}$ goes outside of $K$
is almost surely finite.
\end{theorem}

In order to establish the asymptotic normality of $\wh{\theta}_{n}$,
it is necessary to
introduce a second auxiliary function $\varphi$ defined, for all $t\in
\dR$, by
%
%e2.6 #&#
%
\begin{eqnarray}
\label{defvarphi}
\varphi(t)&=&\dE\biggl[\frac{\sin^2(2\pi(X-t))
}{g^2(X)}\bigl(f^2(X-\theta) + \sigma^2\bigr)\biggr]
\nonumber\\[-8pt]\\[-8pt]
&=& \int_{-1/2}^{1/2}\frac{\sin^2(2\pi(x-t))}{g(x)}\bigl(f^2(x-\theta)+
\sigma^2\bigr) \,dx. \nonumber
\end{eqnarray}
As soon as $4\pi|f_{1}|>1$, denote
%
%e2.7 #&#
%
\begin{equation}
\label{varasympt}
\xi^{2}(\theta)=\frac{\varphi(\theta)}{4\pi|f_{1}|-1}.
\end{equation}
%
%th2.2 #&#
%
\begin{theorem}
\label{thmcltrm}
Assume that~(\ref{equH1}) and~(\ref{equH2}) hold and that
$|\theta|<1/4$. In addition, suppose that $(\veps_{n})$ has a finite
moment of order $>2$
and that $4\pi|f_{1}|>1$. Then, we have the asymptotic normality
%
%e2.8 #&#
%
\begin{equation}
\label{cltrm}
\sqrt{n}(\wh{\theta}_{n}-\theta) \liml\cN(0, \xi^2(\theta)).
\end{equation}
\end{theorem}
%
%re2.1 #&#
%
\begin{rem}
\label{remefficiency}
We clearly have
$
\phi^{\prime}(t)=-2 \pi f_1 \cos(2\pi(\theta-t)).
$
Consequently, the value $\phi^{\prime}(\theta)=-2\pi f_1$ does not
depend upon the unknown parameter $\theta$.
On the one hand, if the first Fourier coefficient $f_1$ of~$f$ is known,
it is possible to provide, via a slight modification of~(\ref{RMA}),
an asymptotically efficient estimator $\wh{\theta}_{n}$ of $\theta$.
More precisely, it is only necessary to replace $\gamma_n=1/n$ in~(\ref{RMA}) by $\gamma_n=\gamma/n$ where
\[
\gamma= \frac{1}{2 \pi| f_1|}.
\]
Then, we deduce from the original work of Fabian~\cite{MR0381189} that
$\wh{\theta}_{n}$ is an asymptotically efficient estimator of $\theta$
with
%
%e2.9 #&#
%
\begin{equation}
\label{cltrmefficient}
\sqrt{n}(\wh{\theta}_{n}-\theta) \liml\cN\biggl(0, \frac{\varphi
(\theta)}{4 \pi^2 f_1^2}\biggr).
\end{equation}
On the other\vspace*{1pt} hand, if $f_1$ is unknown, it is also possible to provide
by the same procedure
an asymptotically efficient estimator $\wh{\theta}_{n}$ of $\theta$
replacing $f_1$ by its natural estimate
\[
\wh{f}_{1,n}= \frac{1}{n} \sum_{k=1}^n \frac{Y_k \cos(2 \pi(X_k -
\wh{\theta}_{k-1}))}{g(X_k)}.
\]
\end{rem}
%
%re2.2 #&#
%
\begin{rem} In the particular case where $4\pi|f_{1}|=1$, it is also
possible to show~\cite{MR1485774} that
\[
\sqrt{\frac{n}{\log(n)}}(\wh{\theta}_{n}-\theta) \liml\cN(0,
\varphi(\theta)).
\]
Asymptotic results are also available when $0<4\pi|f_{1}|<1$.
However, we have chosen to focus our attention on the more attractive case
$4\pi|f_{1}|>1$.
\end{rem}
%
%th2.3 #&#
%
\begin{theorem}
\label{thmlilqsl}
Assume that~(\ref{equH1}) and~(\ref{equH2}) hold and that
$|\theta|<1/4$. In addition, suppose that $(\veps_{n})$ has a finite
moment of order $>2$
and that $4\pi|f_{1}|>1$. Then,
we have the law of iterated logarithm
%
%e2.10 #&#
%
\begin{eqnarray}
\label{lilrm}
\limsup_{n \rightarrow\infty} \biggl(\frac{n}{2 \log\log n} \biggr)^{1/2}
(\wh{\theta}_{n}-\theta)
&=& - \liminf_{n \rightarrow\infty}
\biggl(\frac{n}{2 \log\log n}\biggr)^{1/2}
(\wh{\theta}_{n}-\theta)\nonumber\hspace*{-35pt}\\[-8pt]\\[-8pt]
&=& \xi(\theta)
\qquad\mbox{a.s.}\nonumber\hspace*{-35pt}
\end{eqnarray}
In particular,
%
%e2.11 #&#
%
\begin{equation}
\label{lilsuprm}
\limsup_{n \rightarrow\infty} \biggl(\frac{n}{2 \log\log n} \biggr)
(\wh{\theta}_{n}-\theta)^2=\xi^2(\theta)\qquad
\mbox{a.s.}
\end{equation}
In addition, we also have the quadratic strong law
%
%e2.12 #&#
%
\begin{equation}
\label{qslrm}
\lim_{n\rightarrow\infty}\frac{1}{\log n}\sum_{k=1}^{n}(\wh
{\theta}_{k}-\theta)^{2}=\xi^2(\theta)\qquad
\mbox{a.s.}
\end{equation}
\end{theorem}
\begin{pf}
The proofs are given in Section~\ref{sec5}.
\end{pf}
%
%re2.3 #&#
%
\begin{rem}
\label{remsym}
It is also possible to get rid of the symmetry assumption on~$f$.
However, it requires the knowledge of the first Fourier coefficients of~$f$:
\[
f_1=\int_{-1/2}^{1/2}\cos(2\pi x)f(x) \,dx
\quad\mbox{and}\quad
g_1=\int_{-1/2}^{1/2}\sin(2\pi x)f(x) \,dx.
\]
On the one hand, it is necessary to assume that $f_1\neq0$ or $g_1\neq
0$, and to replace the
first auxiliary function $\phi$ defined in~(\ref{defphi}) by
\begin{eqnarray*}
\Phi(t)&=&f_1\dE\biggl[\frac{\sin(2\pi(X-t)) }{g(X)}f(X-\theta
)\biggr] -g_1\dE\biggl[\frac{\cos(2\pi(X-t)) }{g(X)}f(X-\theta
)\biggr] \\
&=&(f_1^2 + g_1^2)\sin\bigl(2\pi(\theta- t)\bigr).
\end{eqnarray*}
Then, Theorem~\ref{thmascvg} is true for the projected Robbins--Monro algorithm
\[
\wh{\theta}_{n+1}=\pi_{K}(\wh{\theta}_{n}+\gamma
_{n+1}T_{n+1}),
\]
where the initial value $\wh{\theta}_{0} \in K$ and the random
variable $T_{n+1}$ is defined by
\[
T_{n+1}=\frac{f_1\sin(2\pi(X_{n+1}-\wh{\theta
}_{n}))}{g(X_{n+1})}Y_{n+1}-\frac{g_1\cos(2\pi(X_{n+1}-\wh{\theta
}_{n}))}{g(X_{n+1})}Y_{n+1}.
\]
On the other hand, we also have to replace the second function $\varphi
$ defined in~(\ref{defvarphi}) by
\begin{eqnarray*}
\Psi(t)&=&\dE\biggl[\frac{(f_1\sin(2\pi(X-t)) -g_1\cos(2\pi
(X-t)))^2}{g^2(X)}\bigl(f^2(X-\theta) + \sigma^2\bigr)\biggr] \\
&=& \int_{-1/2}^{1/2}\frac{(f_1\sin(2\pi(x-t)) -g_1\cos(2\pi
(x-t)))^2}{g(x)}\bigl(f^2(x-\theta)+ \sigma^2\bigr) \,dx.
\end{eqnarray*}
Then, as soon as $4\pi(f_{1}^2+g_{1}^2)>1$, Theorems~\ref{thmcltrm}
and~\ref{thmlilqsl} hold with
\[
\xi^{2}(\theta)=\frac{\Psi(\theta)}{4\pi(f_{1}^2+g_{1}^2) -1}.
\]
In the rest of the paper, we shall not go in that direction as our
strategy is to make very few assumptions on the Fourier coefficients
of~$f$.
\end{rem}
%

%%%%%%%%%%%%%%%%%%%%%%%%%%%%%%%%%%%%%%%%%%%%%%%%%%%%%%%%%%%%%%%%%%%%%%%%%%%%%%%%%%%%

%s3 #&#
\section{Estimation of the regression function}\label{sec3}

%%%%%%%%%%%%%%%%%%%%%%%%%%%%%%%%%%%%%%%%%%%%%%%%%%%%%%%%%%%%%%%%%%%%%%%%%%%%%%%%%%%%
This section is devoted to the nonparametric estimation of the
regression function~$f$ via a recursive Nada\-raya--Watson
estimator. On the one hand, we add the standard hypothesis:
{\renewcommand{\theequation}{$\mathcal{H}_{3}$}
\begin{equation}\label{equH3}
\mbox{The regression function~$f$ is Lipschitz}.\hspace*{110pt}
\end{equation}}

\vspace*{-\baselineskip}

\noindent
On the other hand, we recall that under~(\ref{equH2}), the
function~$f$ is assumed to be symmetric.
Consequently, we follow the same approach as the one developed by Stone
\cite{MR0362669} for the estimation of
a symmetric probability density function replacing the estimator~(\ref{RNW})
by its symmetrized version
%
%e3.1 #&#
%
\setcounter{equation}{0}
\begin{equation}
\label{RNWS}
\wh{f}_{n}(x)=\frac{\sum_{k=1}^{n} (W_{k}(x)+W_{k}(-x))Y_{k}}{\sum
_{k=1}^{n} (W_{k}(x)+W_{k}(-x))},
\end{equation}
where
\[
W_{n}(x)=\frac{1}{h_{n}}K\biggl(\frac{X_{n}-\widehat{\theta
}_{n-1}-x}{h_{n}}\biggr).
\]
The bandwidth $(h_n)$ is a sequence of positive real numbers,
decreasing to zero, such that $n h_n$ tends to infinity. For the sake
of simplicity, we propose to make use of
$h_n = 1/n^{\alpha}$ with $\alpha\in\ ]0,1[$.
Moreover, we shall assume in all the sequel that the kernel $K$ is a
positive symmetric function, bounded with compact support,
twice differentiable with bounded derivatives, satisfying
\[
\int_{\dR} K(x) \,dx = 1
\quad\mbox{and}\quad
\int_{\dR} K^2(x) \,dx= \nu^2.
\]

Our next result deals with the almost sure convergence of the estimator~$\wh{f}_{n}$.
%
%th3.1 #&#
%
\begin{theorem}
\label{thmaspnw}
Assume that~(\ref{equH1}),~(\ref{equH2}) and
(\ref{equH3}) hold and that $|\theta|<1/4$
and the sequence $(\veps_{n})$ has a finite moment of order $>2$.
Then, for any $x\in\dR$ such that $|x|\leq{1/2}$,
%
%e3.2 #&#
%
\begin{equation}
\label{Cvgfchapn}
\lim_{n\rightarrow\infty}
\wh{f}_{n}(x)=f(x) \qquad\mbox{a.s.}
\end{equation}
\end{theorem}

The asymptotic normality of the estimator $\wh{f}_{n}$ is as follows.
%
%th3.2 #&#
%
\begin{theorem}
\label{thmcltnw}
Assume that~(\ref{equH1}),~(\ref{equH2}) and~(\ref{equH3})
hold and that $|\theta|<1/4$
and the sequence $(\veps_{n})$ has a finite moment of order $>2$.
Then, as soon as
the bandwidth $(h_n)$ satisfies $h_n = 1/n^{\alpha}$ with $\alpha
>1/3$, we have
for any $x\in\dR$ such that
$|x|\leq{1/2}$ with $x\neq{0}$, the pointwise asymptotic normality
%
%e3.3 #&#
%
\begin{equation}
\label{cltnwx}
\sqrt{nh_n}\bigl(\wh{f}_{n}(x)-f(x)\bigr) \liml\cN\biggl(0,\frac{\sigma
^{2}\nu^2}{(1+\alpha)(g(\theta+x)+g(\theta-x))}\biggr).
\end{equation}
In addition, for $x=0$,
%
%e3.4 #&#
%
\begin{equation}
\label{cltnwzero}
\sqrt{nh_n}\bigl(\wh{f}_{n}(0)-f(0)\bigr) \liml\cN\biggl(0,\frac{\sigma
^{2}\nu^2}{(1+\alpha)g(\theta)}\biggr).
\end{equation}
\end{theorem}
\begin{pf}
The proofs are given in Section~\ref{sec6}.
\end{pf}

%%%%%%%%%%%%%%%%%%%%%%%%%%%%%%%%%%%%%%%%%%%%%%%%%%%%%%%%%%%%%%%%%%%%%%%%%%%%%%%%%%%%

%s4 #&#
\section{Simulations}\label{sec4}

%%%%%%%%%%%%%%%%%%%%%%%%%%%%%%%%%%%%%%%%%%%%%%%%%%%%%%%%%%%%%%%%%%%%%%%%%%%%%%%%%%%%

The goal of this section is to illustrate via some numerical
experiments the good performances of our estimation
strategy. The first subsection is devoted to simulated data created
according to the model~(\ref{Sempar})
while the second one deals with real ECG data taken from the MIT-BIH
database. Our aim is to propose an
efficient and easy to handle procedure in order to detect atrial
fibrillation using ECG records.
An interesting study on ECG analysis in order to detect cardiac
arrhythmia may also be found in~\cite{Trigano}.

%%%%%%%%%%%%%%%%%%%%%%%%%%%%%%%%%%%%%%%%%%%%%%%%%%%%%%%%%%%%%%%%%%%%%%%%%%%%%%%%%%%%

%s4.1 #&#
\subsection{Simulated data}\label{sec41}

%%%%%%%%%%%%%%%%%%%%%%%%%%%%%%%%%%%%%%%%%%%%%%%%%%%%%%%%%%%%%%%%%%%%%%%%%%%%%%%%%%%%

\begin{figure}

\includegraphics{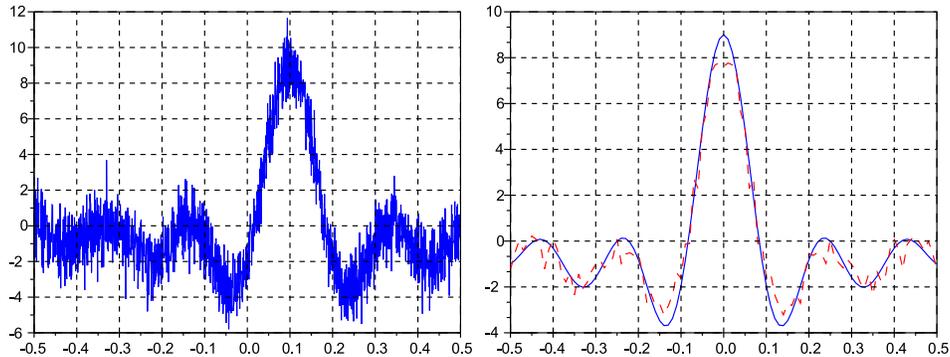}

\caption{Simulated data and almost sure convergence.}
\label{data}
\end{figure}

Consider the semiparametric regression model
\[
Y_{n}=f(X_{n}-\theta)+\veps_n,
\]
where $\theta=1/10$ and the periodic shape function~$f$ is given, for
$p\geq1$ and for all $x \in\dR$, by
\[
f(x)=\sum_{k=1}^{p}\cos(2k\pi x)
\]
with $p=8$. We have chosen $(X_n)$ and $(\veps_n)$ as two independent
sequences of independent random variables
with $\mathcal{U}[-1/2,1/2]$ and $\cN(0,1)$ distributions, respectively.
The simulated data are given in the left-hand side of Figure~\ref{data}.

%
%f1 #&#
%

\begin{figure}

\includegraphics{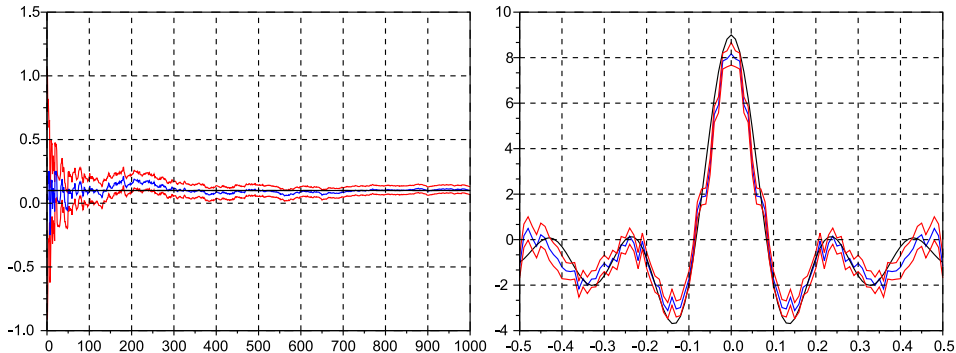}

\caption{Confidence intervals for $\theta$ and~$f$.}
\label{IC}
\end{figure}

For the estimation of the shift parameter $\theta$, we implement our
Robbins--Monro procedure
with $n=1000$ iterations. We obtain the estimate $\wh{\theta}_{n}=
0.1014$ which shows the good asymptotic behavior
of the estimator $\wh{\theta}_{n}$ comparing to the true value
$\theta=1/10$. Moreover, using convergence~(\ref{cltrm}), one can obtain
confidence intervals for the shift parameter. More precisely, they are
given, for all $n\geq1$, by
\[
I_{n}(\theta)=\biggl[\wh{\theta}_{n}-q_{\beta}\frac{\widehat{\xi
}_{n}(\theta)}{\sqrt{n}},\wh{\theta}_{n}+q_{\beta}\frac{\widehat
{\xi}_{n}(\theta)}{\sqrt{n}}\biggr],
\]
where $q_{\beta}$ stands for the quantile of order $0<\beta<1$ of the
$\cN(0,1)$ distribution and $\widehat{\xi}_{n}(\theta)$ is a
consistent estimator of $\xi(\theta)$
given by~(\ref{varasympt}). In our particular case, it is not
necessary to estimate $\xi^{2}(\theta)$ since via straightforward
calculations, $f_{1}=1/2$ and
\[
\xi^{2}(\theta)=\frac{7}{8(2\pi-1)}.
\]
Moreover, for $n=1000$ and for a risk $\beta=5\%$, the confidence
interval is precisely $I_{n}(\theta)=[0.0762;0.1266]$. The length of
$I_{n}(\theta)$ is 0.0504, which is rather small, so our Robbins--Monro
procedure performs pretty well. All confidence intervals $I_{n}(\theta
)$, for $n=1,\ldots,1000$, are drawn in red in the left-hand side of Figure
\ref{IC}.

For the estimation of the regression function~$f$,
we make use of the uniform kernel $K$ on the interval $[-1,1]$, and the
bandwidth $h_n=1/n^{\alpha}$ with $\alpha=9/10$.
In addition, it follows from convergences~(\ref{cltnwx}) and~(\ref{cltnwzero}) that for $n=1000$ and for all $x\in{[-1/2,1/2]}$, a
confidence interval for $f(x)$ is given by
\[
J_{n}(x)=\biggl[\wh{f}_{n}(x)-q_{\beta}\frac{\widehat{v}_{n}(x,\wh
{\theta}_{n})}{\sqrt{n h_{n}}},\wh{f}_{n}(x)+q_{\beta}\frac
{\widehat{v}_{n}(x,\wh{\theta}_{n})}{\sqrt{n h_{n}}}\biggr],
\]
where $q_{\beta}$ stands for the quantile of order $0<\beta<1$ of the
$\cN(0,1)$ distribution and $\widehat{v}_{n}^{ 2}(x,\wh{\theta
}_{n})$ is a consistent estimator of the asymptotic variance
$v^{2}(x,\theta)$ in Theorem~\ref{thmcltnw}.
In our particular case, $\nu^{2}=1/2$ and
\[
v^{2}(x,\theta) = \cases{
5/19, &\quad if $-1/2\leq{x}< -2/$ or $2/5<x\leq1/2$, \cr
5/38, &\quad if $-2/5\leq{x}\leq2/5$ and $x\neq0$,\cr
5/19, &\quad if $x=0$.}
\]

All confidence intervals $J_{n}(x)$, for all $x\in{}[-1/2,1/2]$, are
drawn in red in the right-hand side of Figure~\ref{IC}.
On the one hand, the simulations show that the largest length of the
confidence intervals $J_{n}(x)$
is for $x=-0.47$ and $x=0.47$ and the length is precisely equal to
$1.0066$. On the other hand, the smallest length of the confidence
intervals $J_{n}(x)$ is for $x=-0.04$ and $x=0.04$ and is equal to
$0.7118$. The fact that there are two values of $x$ for the largest and
the smallest length of confidence intervals is due to the symmetry of
the estimator $\wh{f}_{n}$.
Then, one can observe on this first set of simulated data that the
Robbins--Monro estimator $\wh{\theta}_{n}$ of $\theta$ as well as
the Nadaraya--Watson estimator $\wh{f}_{n}$ of~$f$ perform pretty well.

%
%f2 #&#
%

Our second experiment deals with $30$ curves according to the model
\[
Y_{n}=f(X_{n}-\theta)+\veps_n,
\]
where $\theta=-1/5$ for the first $10$ curves and $\theta^{\prime
}=1/10$ for the last $20$ curves.
The periodic shape function~$f$ is given, for all $x \in[-1/2,1/2]$, by
\[
f(x)=\cos(2\pi x)+\sin(2\pi x)+\cos(2\pi x)\sin(2\pi x).
\]

\begin{figure}

\includegraphics{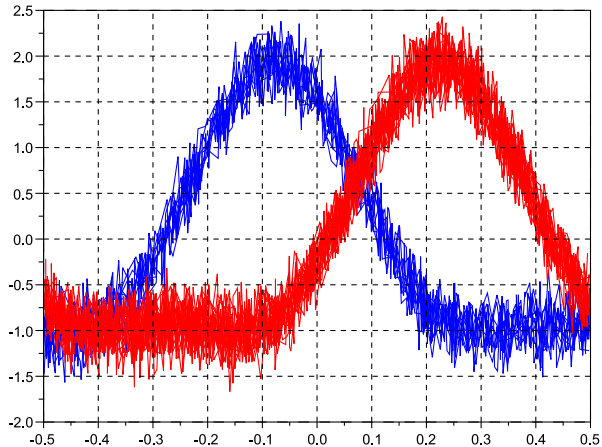}

\caption{Simulated data with two different values $\theta$ and
$\theta^{\prime}$.}
\label{data2}
\end{figure}

Our goal is to propose a statistical procedure in order to detect a lag
between the first $10$ curves with $\theta=-1/5$ and
the last $20$ curves with $\theta^{\prime}=1/10$. In other words, we
want to observe whether or not the value
$\Delta= \theta^{\prime}- \theta$ is far away from zero.
We have chosen $(X_n)$ and $(\veps_n)$ as two independent sequences of
independent and Gaussian random variables
with uniform distribution on $[-1/2,1/2]$ and $\cN(0,1/5)$
distribution, respectively. Each curve is drawn with $n=200$ points.
The different curves are given in Figure~\ref{data2}.\vadjust{\goodbreak}

%
%f3 #&#
%

On the one hand, we estimate the first value $\theta=-1/5$ from the
first $10$ curves.
We implement our Robbins--Monro procedure with $n=200$ iterations for
the first estimate
$\wh{\theta}_{n}$ of $\theta$ evaluated on the first curve, then
with $n=400$ iterations for the second
estimate $\wh{\theta}_{n}$ of $\theta$ evaluated on the two first
curves, and so on, until the calculation
of the last estimate $\wh{\theta}_{n}$ of $\theta$ with $n=2000$.
Therefore, we obtain $-0.1950$ for the
arithmetic mean of the first $10$ estimates $\wh{\theta}_{n}$ of
$\theta$. We continue with the same
procedure on all the set of curves. The value of the eleven estimates
with $n=2200$ is $0.0986$.
This value is significantly different from the first $10$ estimates. It
corresponds to the first curve
simulated with $\theta^{\prime}=1/10$. Furthermore, we obtain
$0.0998$ for the
arithmetic mean of the last $20$ estimates $\wh{\theta}_{n}$ of
$\theta^{\prime}$.
Finally, our statistical procedure allows us to detect a change of
parameterization
from the value $\theta=-1/5$ to the value $\theta^{\prime}=1/10$ as
$\wh{\Delta}_n= 0.0998+0.1950=0.2948$. In order to compute more
accurate values
of $\wh{\theta}_{n}$, one can replace $\gamma_n=1/n$ in~(\ref{RMA})
by $\gamma_n=1/n^{a}$
where $1/2 < a <1$. This will be done for the implementation of our
Robbins--Monro procedure
on real ECG data.

%%%%%%%%%%%%%%%%%%%%%%%%%%%%%%%%%%%%%%%%%%%%%%%%%%%%%%%%%%%%%%%%%%%%%%%%%%%%%%%%%%%%

%s4.2 #&#
\subsection{Real ECG data}\label{sec42}

%%%%%%%%%%%%%%%%%%%%%%%%%%%%%%%%%%%%%%%%%%%%%%%%%%%%%%%%%%%%%%%%%%%%%%%%%%%%%%%%%%%%

\begin{figure}

\includegraphics{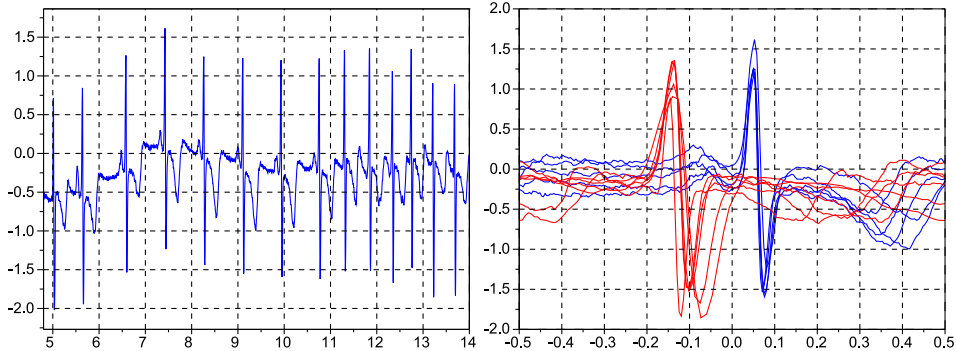}

\caption{Original data.}
\label{signal}
\vspace*{-3pt}
\end{figure}

We shall now focus our attention on real ECG data. We have chosen the
record $04015$ in the Atrial Fibrillation (AF) database provided by
MIT-BIH database.
Each recording consists in a continuous digitized ECG signal measured
over 1 hour in order to detect AF which is the most common
cardiac arrythmia. A stronger indicator of AF is the absence of P waves
or the irregularities of RR interval on an electrocardiogram.
We refer the reader to~\cite{Clifford} for an interesting book on
statistical methods and tools for ECG data analysis.
Our aim is to propose a statistical procedure in order to detect
irregularities of RR interval on the ECG record $04015$.
The record and its\vadjust{\goodbreak} projection on the interval $[-1/2,1/2]$ are given in
Figure~\ref{signal}. The size of the data set is $2038$.
We assume that the model
\[
Y_{n}=f(X_{n}-\theta)+\veps_n
\]
fits the data, where the sequence $(X_n)$ is uniformly distributed over
the interval $[-1/2,1/2]$.
The periodic shape function~$f$ is clearly not symmetric. However, we
already saw in Remark~\ref{remsym} that
our Robbins--Monro procedure still holds for nonsymmetric regression function.

%
%f4 #&#
%

%
%f5 #&#
%
\begin{figure}[b]

\includegraphics{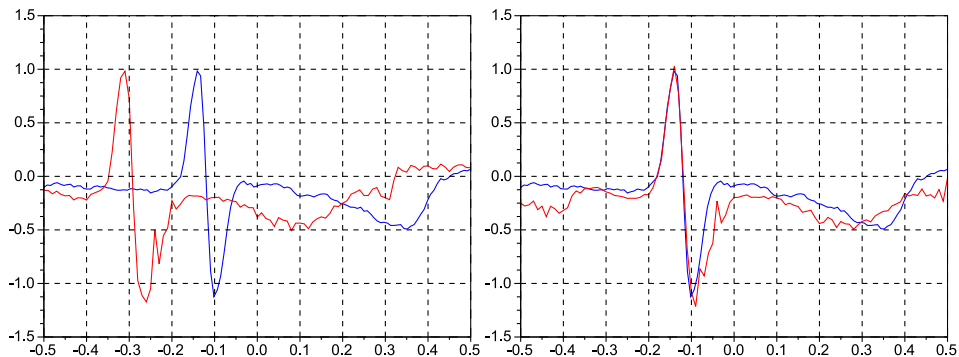}

\caption{Reconstruction of the ECG.}
\label{reconstruction}
\end{figure}

As for simulated data, in view of the signal, we would find two
different values $\theta$ and $\theta^{\prime}$.
The first value $\theta$ is associated with the first part of the
signal, while the second value $\theta^{\prime}$
corresponds to the second part. The difference $\Delta= \theta
^{\prime}- \theta$
between the two parameters would explain the lag between the two parts
of the signal.
A~value of $\Delta$ far away from zero could be interpreted as the
detection of
irregularities of RR interval which confirms the diagnostic of atrial
fibrillation.
On this record, our Robbins--Monro procedure with $n=800$ iterations
leads to the first estimate $\wh{\theta}_{n}=0.1734$ for $\theta$
and the last estimate $\wh{\theta}_{n}=-0.0092$ for $\theta^{\prime
}$ with
$n=1238$. The value $\wh{\Delta}_n= -0.0092-0.1734=-0.1826$
explains the lag in Figure~\ref{signal}.\vadjust{\goodbreak}
Figure~\ref{reconstruction} shows that our Nadaraya--Watson procedure for the
reconstruction of ECG signals works pretty well.

%%%%%%%%%%%%%%%%%%%%%%%%%%%%%%%%%%%%%%%%%%%%%%%%%%%%%%%%%%%%%%%%%%%%%%%%%%%%%%%%%%%%

%s5 #&#
\section{Proofs of the parametric results}\label{sec5}

%%%%%%%%%%%%%%%%%%%%%%%%%%%%%%%%%%%%%%%%%%%%%%%%%%%%%%%%%%%%%%%%%%%%%%%%%%%%%%%%%%%%

%s5.1 #&#
\subsection{\texorpdfstring{Proof of Theorem \protect\ref{thmascvg}}{Proof of Theorem 2.1}}\label{sec51}
We can assume without loss of generality that \mbox{$f_1>0$} inasmuch as the proof
for $f_1<0$ follows exactly the same lines.
Denote by $\cF_{n}$ the $\sigma$-algebra of the events occurring up
to time $n$,
$\cF_{n}=\sigma(X_0,\veps_0,\ldots, X_n,\veps_n)$. First of all,
we shall calculate the two first conditional
moments of the random variable $T_n$ given by~(\ref{DefT}). It follows
from~(\ref{Sempar}) that
\begin{eqnarray*}
\dE[T_{n+1}|\mathcal{F}_{n}]&=&\dE\biggl[\frac{\sin(2\pi
(X_{n+1}-\wh{\theta}_{n}))Y_{n+1}}{g(X_{n+1})}\Big|\mathcal{F}_{n}
\biggr]\\
&=&\dE\biggl[\frac{\sin(2\pi(X_{n+1}-\wh{\theta
}_{n}))(f(X_{n+1}-\theta)+\veps_{n+1})}{g(X_{n+1})}\Big|\mathcal
{F}_{n}\biggr].
\end{eqnarray*}
On the one hand, as $(X_n)$ is a sequence of independent random
variables sharing the same distribution as a random
variable $X$, we have
%
%e5.1 #&#
%
\begin{equation}
\label{meTn}\qquad
\dE\biggl[\frac{\sin(2\pi(X_{n+1}-\wh{\theta
}_{n}))f(X_{n+1}-\theta)}{g(X_{n+1})}\Big|\mathcal{F}_{n}\biggr]=\phi
(\wh{\theta}_{n})\qquad
\mbox{a.s.},
\end{equation}
where $\phi$ is the function given by~(\ref{defphi}). On the other
hand, as $(X_n)$ and $(\veps_n)$ are two independent sequences
and $(\veps_n)$ is a sequence of independent and square integrable
random variables with zero mean, we also have
\[
\dE\biggl[\frac{\sin(2\pi(X_{n+1}-\wh{\theta}_{n}))\veps
_{n+1}}{g(X_{n+1})}\Big|\mathcal{F}_{n}\biggr]=\dE\biggl[\frac{\sin
(2\pi(X-\wh{\theta}_{n}))}{g(X)}\biggr]
\dE[\veps_{n+1}]=0.
\]
Hence,~(\ref{meTn}) leads to
%
%e5.2 #&#
%
\begin{equation}
\label{meanTn}
\dE[T_{n+1}|\mathcal{F}_{n}]=\phi(\wh{\theta}_{n})\qquad
\mbox{a.s.}
\end{equation}
On the other hand,
\begin{eqnarray*}
T_{n+1}^2&=&\frac{\sin^{2}(2\pi(X_{n+1}-\wh{\theta
}_{n}))Y_{n+1}^2}{g^{2}(X_{n+1})}\\
&=&\frac{\sin^{2}(2\pi(X_{n+1}-\wh{\theta
}_{n}))(f^{2}(X_{n+1}-\theta)+2\veps_{n+1}f(X_{n+1}-\theta)+\veps
_{n+1}^2 )}{g^{2}(X_{n+1})}.
\end{eqnarray*}
Consequently, as the function~$f$ is bounded, the density $g$ is
positive on $[-1/2,1/2]$, and
$\dE[\veps_{n+1}^2|\mathcal{F}_{n}]=\dE[\veps_{n+1}^2]=\sigma^2$,
we obtain that
%
%e5.3 #&#
%
\begin{equation}
\label{msTn}\qquad
\dE[T_{n+1}^2|\mathcal{F}_{n}]=\dE\biggl[\frac{\sin^2(2\pi(X-\wh
{\theta}_{n})) }{g^2(X)}\bigl(f^2(X-\theta) + \sigma^2\bigr)\biggr]=\varphi
(\wh{\theta}_{n}),
\end{equation}
where $\varphi$ is given by~(\ref{defvarphi}).
Therefore, as~$f$ is bounded and $g$ does not vanish on its support
$[-1/2,1/2]$, we deduce from~(\ref{msTn}) that
for some constant $M>0$
%
%e5.4 #&#
%
\begin{equation}
\label{meansquareTn}
\sup_{n\geq0}
\dE[T_{n+1}^2|\mathcal{F}_{n}]\leq M
\qquad\mbox{a.s.}
\end{equation}
Furthermore, for all $n\geq0$, let $V_{n}=(\wh{\theta}_{n}-\theta
)^2$. We clearly have
\begin{eqnarray*}
V_{n+1}&=&(\wh{\theta}_{n+1}-\theta)^2\\
&=&\bigl(\pi_{K}(\wh{\theta}_{n}+\gamma_{n+1}T_{n+1})-\theta\bigr)^2\\
&=&\bigl(\pi_{K}(\wh{\theta}_{n}+\gamma_{n+1}T_{n+1})-\pi_{K}(\theta)\bigr)^2
\end{eqnarray*}
as we have assumed that $\theta$ belongs to $K$. Since $\pi_{K}$ is a
Lipschitz function with Lipschitz constant $1$, we obtain that
\begin{eqnarray*}
V_{n+1}
&\leq&(\wh{\theta}_{n}+\gamma_{n+1}T_{n+1}-\theta)^2\\
&\leq&V_{n}+\gamma_{n+1}^{2}T_{n+1}^{2}+2\gamma_{n+1}T_{n+1}(\wh
{\theta}_{n}-\theta).
\end{eqnarray*}
Hence, it follows from~(\ref{meanTn}) and~(\ref{meansquareTn}) that
%
%e5.5 #&#
%
\begin{eqnarray}\label{RSiegmund}
\dE[V_{n+1}|\mathcal{F}_{n}]&\leq&V_{n}+\gamma_{n+1}^{2}\dE
[T_{n+1}^{2}|\mathcal{F}_{n}]+2\gamma_{n+1}(\wh{\theta}_{n}-\theta
)\dE[T_{n+1}|\mathcal{F}_{n}] \nonumber\hspace*{-30pt}\\[-8pt]\\[-8pt]
&\leq& V_{n}+\gamma_{n+1}^{2}M+2\gamma_{n+1}(\wh{\theta
}_{n}-\theta)\phi(\wh{\theta}_{n})\qquad
\mbox{a.s.}
\nonumber\hspace*{-30pt}
\end{eqnarray}
In addition, as $\wh{\theta}_{n} \in K$, $|\wh{\theta}_{n}|<1/4$,
$|\wh{\theta}_{n}-\theta|<1/2$ which implies that $(\wh{\theta
}_{n}-\theta)\phi(\wh{\theta}_{n})<0$.
Then, we deduce from~(\ref{RSiegmund}) together with the
Robbins--Siegmund theorem
(see Duflo~\cite{MR1485774}, page 18) that the sequence $(V_n)$
converges a.s. to a~finite random variable $V$
and
%
%e5.6 #&#
%
\begin{equation}
\label{sumRS}
\sum_{n=1}^\infty\gamma_{n+1}(\theta- \wh{\theta}_{n})\phi(\wh
{\theta}_{n})<+\infty\qquad\mbox{a.s.}
\end{equation}
Assume by contradiction that $V \neq0$ a.s. Then, one can find
$0<a<b<1/2$ such that, for
$n$ large enough, the event $\{a<| \wh{\theta}_{n}-\theta|<b\}$ is
not negligible. However, on this annulus, one can also find some
constant $c>0$ such that $(\theta- \wh{\theta}_{n})\phi(\wh{\theta
}_{n})>c$ which, by~(\ref{sumRS}), implies that
\[
\sum_{n=1}^\infty\gamma_{n}<+\infty.
\]
This is of course in contradiction with assumption~(\ref{hypgamma}).
Consequently, it follows that $V=0$ a.s.
leading to the almost sure convergence of $\wh{\theta}_{n}$ to
$\theta$.
\newline
It remains to show that
$\wh{\theta}_{n}+\gamma_{n+1}T_{n+1}$ goes almost surely outside of
$K$ a~finite number of times.
For all $n\geq1$, denote
\[
N_n=\sum_{k=0}^{n-1} \rI_{\{| \wh{\theta}_{k}+\gamma
_{k+1}T_{k+1}|>1/4\}}.\vadjust{\goodbreak}
\]
The random sequence $(N_n)$ is nondecreasing. Assume by
contradiction\break
that~$N_n$ goes to infinity a.s.
Then, one can find a subsequence $(n_k)$ such that $(N_{n_k})$ is
increasing. Consequently, for all $n_k>0$,
\[
|\wh{\theta}_{n_k}+\gamma_{n_k+1}T_{n_k+1}|>\tfrac{1}{4}
\qquad\mbox{a.s.},
\]
which implies that $ |\wh{\theta}_{n_k +1}|=1/4$ a.s. Hence,
\[
\lim_{n_k\rightarrow\infty} |\wh{\theta}_{n_k}|= |\theta|=\frac{1}{4}
\qquad\mbox{a.s.}
\]
leading to a contradiction as $|\theta|<1/4$. Finally, $(N_n)$
converges to a finite limiting value a.s.
which completes the proof of Theorem~\ref{thmascvg}.

%s5.2 #&#
\subsection{\texorpdfstring{Proof of Theorem \protect\ref{thmcltrm}}{Proof of Theorem 2.2}}\label{sec52}

We assume without loss of generality that \mbox{$f_{1}>0$}.
Our goal is to apply Theorem 2.1 of Kushner and Yin
(\cite{MR1993642}, page~330). First of all, as $\gamma_{n}=1/n$,
the condition on the decreasing step is satisfied. Moreover, we already
saw that $\wh{\theta}_n$ converges almost surely to $\theta$.
Consequently, all the local assumptions of Theorem 2.1 of
\cite{MR1993642} are satisfied.
In addition, it follows from $(\ref{meanTn})$ that $\dE
[T_{n+1}|\cF_{n}]=\phi(\wh{\theta}_n)$ a.s. and the function~$\phi$ is
continuously differentiable since $\phi(t)=f_{1}\sin(2\pi(\theta-t))$.
Hence, $\phi(\theta)=0$ and $\phi^{\prime}(\theta)=-2\pi f_{1}$
and $4\pi f_{1}>1$.
Furthermore, we deduce from~$(\ref{msTn})$ that
\[
\dE[T_{n+1}^{2}|\cF_{n}]=\varphi(\wh{\theta
}_n) \qquad\mbox{a.s.},
\]
which leads to
\[
\lim_{n\rightarrow\infty}\dE[T_{n+1}^{2}|\cF_{n}
]=\varphi(\theta) \qquad\mbox{a.s.}
\]
Consequently, if we are able to prove that the sequence $(W_n)$ given by
\[
W_{n}=\frac{(\wh{\theta}_{n}-\theta)^{2}}{\gamma_{n}}
\]
is tight, then we shall deduce from Theorem 2.1 of~\cite{MR1993642} that
\[
\sqrt{n}(\wh{\theta}_{n}-\theta) \liml\cN(0, \xi^2(\theta)),
\]
where
\[
\xi^2(\theta)=\varphi(\theta)\int_{0}^{+\infty}\exp\bigl((1-4\pi
f_{1})t\bigr) \,dt=\frac{\varphi(\theta)}{4\pi f_{1}-1}.
\]
Therefore, it remains to prove the tightness of the sequence $(W_{n})$.
It follows from~(\ref{RSiegmund}) that for some constant $M>0$ and for
all $n\geq1$,
%
%e5.7 #&#
%
\begin{equation}
\label{delta}
\mathbb{E}[W_{n+1}|\mathcal{F}_{n}]\leq{(1+\gamma_{n})W_{n}+\gamma
_{n+1}M+2(\wh{\theta}_{n}-\theta)\phi(\wh{\theta}_{n})}.
\end{equation}
Moreover, we have for all $x\in\dR$, $\phi(x)=2\pi f_{1}(\theta
-x)+f_{1}(\theta-x)v(x)$ where
\[
v(x)=\frac{\sin(2 \pi(\theta-x) )-2 \pi(\theta-x)}{(\theta-x)}.\vadjust{\goodbreak}
\]
By the continuity of the function $v$, one can find $0<\varepsilon
<1/2$ such that, if $|x -\theta|< \varepsilon$,
%
%e5.8 #&#
%
\begin{equation}
\label{majv}
\frac{q}{2f_1} < v(x)<0.
\end{equation}
We also deduce from~(\ref{delta}) that for all $n\geq1$,
%
%e5.9 #&#
%
\begin{equation}
\label{delta2}
\mathbb{E}[W_{n+1}|\mathcal{F}_{n}]\leq{W_{n}+2\gamma
_{n}W_{n}\bigl(q-f_{1}v(\wh{\theta}_{n})\bigr)+\gamma_{n}M}
\end{equation}
with $2q=1-4\pi f_{1}$ which means that $q<0$.
Moreover, let $A_n$ and $B_n$ be the sets $A_{n}=\{|\wh{\theta
}_{n}-\theta|\leq\varepsilon\}$
and
\[
B_n={ \bigcap_{k=m}^n} A_k
\]
with $1\leq m\leq n$. Then, it follows from~(\ref{majv}) that
%
%e5.10 #&#
%
\begin{equation}
\label{majvv}
0<-f_1 v(\wh{\theta}_{n})\rI_{B_{n}}< -\biggl(\frac{q}{2}\biggr)\rI
_{B_{n}}.
\end{equation}
Hence, we deduce from the conjunction of~(\ref{delta2})
and~(\ref{majvv}) that for all $n\geq m$,
%
%e5.11 #&#
%
\begin{eqnarray}\label{IF}
\dE[W_{n+1}\rI_{B_{n}}|\cF_{n}]&\leq& W_{n}\rI_{B_{n}} + 2 \gamma
_n W_n \rI_{B_{n}} \biggl( q -\frac{q}{2} \biggr)
+\gamma_{n}M \nonumber\\[-8pt]\\[-8pt]
&\leq& W_{n}\rI_{B_{n}}(1+q\gamma_{n})+\gamma_{n}M.\nonumber
\end{eqnarray}
Since $B_{n+1}=B_n \cap A_{n+1}$, $B_{n+1}\subset{B_{n}}$, and we
obtain by taking the expectation on both sides
of~(\ref{IF}) that for all $n\geq m$,
%
%e5.12 #&#
%
\begin{equation}
\label{II}
\dE[W_{n+1}\rI_{B_{n+1}}]\leq(1+q\gamma_{n})\dE[W_{n}\rI
_{B_{n}}]+\gamma_{n}M.
\end{equation}
From now on, denote $\alpha_{n}=\dE[W_{n}\rI_{B_{n}}]$. We infer
from~(\ref{II}) that for all $n\geq m$,
%
%e5.13 #&#
%
\begin{equation}
\label{alpha}
\alpha_{n+1}\leq\beta_{n}\alpha_{m}+M\beta_{n}\sum_{k=m}^{n}\frac
{\gamma_{k}}{\beta_{k}}
\qquad\mbox{where }
\beta_{n}=\prod_{k=m}^{n}(1+q\gamma_{k}).
\end{equation}
As $\gamma_n=1/n$, it follows from straightforward calculations that
$\beta_{n}=O(n^q)$ and
\[
\sum_{k=1}^{n}\frac{\gamma_{k}}{\beta_{k}}=O(n^{-q}).
\]
Consequently,~(\ref{alpha}) immediately leads to
%
%e5.14 #&#
%
\begin{equation}
\label{majalpha}
\sup_{n\geq m}\alpha_n<+\infty.
\end{equation}
We are now in position to prove the tightness of the sequence $(W_{n})$.
Indeed, it was already proved in
Theorem~\ref{thmascvg} that $\wh{\theta}_{n}$ converges to $\theta$ a.s.
Consequently, if
\[
C_n=\bigcup_{k\geq n} \overline{A}_k,\vadjust{\goodbreak}
\]
then $\dP(C_n)$ converges to zero as $n$ tends to infinity. Moreover,
for $n\geq m$, $\overline{B}_n \subset C_m$ which implies that as
$m,n$ tend to infinity,
$\dP(\overline{B}_n)$ goes to zero. For all $\xi, K >0$ and for all
$n\geq m$ with $m$ large enough,
%
%e5.15 #&#
%
\begin{eqnarray}\label{majPWn}
\dP(W_{n}>K)&\leq& \dP(W_{n}\rI_{B_{n}}>K/2)+\dP(W_{n}\rI
_{\overline{B}_{n}}>K/2) \nonumber\\[-8pt]\\[-8pt]
&\leq&\frac{2}{K}\dE[W_{n}\rI_{B_{n}}]+\dP(\overline{B}_{n}).
\nonumber
\end{eqnarray}
We deduce from~(\ref{majalpha}) that one can find $K$ depending on
$\xi$ such that the first term on the right-hand side of~(\ref{majPWn}) is smaller than
$\xi/2$. It is also the case for the second term as $\dP(\overline
{B}_n)$ goes to zero. Finally,
for all $\xi>0$, it exists $K>0$ such that for $m$ large enough,
\[
\sup_{n \geq m}\dP(W_{n}>K)<\xi,
\]
which implies the tightness of $(W_{n})$ and completes the proof of
Theorem~\ref{thmcltrm}.

%s5.3 #&#
\subsection{\texorpdfstring{Proof of Theorem \protect\ref{thmlilqsl}}{Proof of Theorem 2.3}}\label{sec53}
As the number of times that the random variable $\wh{\theta
}_{n}+\gamma_{n+1}T_{n+1}$ goes outside of $K$
is almost surely finite, the sequence $(\wh{\theta}_n)$ shares the
same almost sure
asymptotic properties as the classical Robbins--Monro algorithm.
Consequently, we deduce the law of iterated logarithm
given by~(\ref{lilrm}) from Theorem 1 of~\cite{MR0365954}; see also
Hall and Heyde (\cite{MR624435}, page~240),
and the quadratic strong law given by~(\ref{qslrm}) from Theorem~3 of
\cite{MR1654569}.

%%%%%%%%%%%%%%%%%%%%%%%%%%%%%%%%%%%%%%%%%%%%%%%%%%%%%%%%%%%%%%%%%%%%%%%%%%%%%%%%%%%%

%s6 #&#
\section{Proofs of the nonparametric results}\label{sec6}

%%%%%%%%%%%%%%%%%%%%%%%%%%%%%%%%%%%%%%%%%%%%%%%%%%%%%%%%%%%%%%%%%%%%%%%%%%%%%%%%%%%%

%s6.1 #&#
\subsection{\texorpdfstring{Proof of Theorem \protect\ref{thmaspnw}}{Proof of Theorem 3.1}}\label{sec61}
In order to prove the almost sure pointwise convergence of Theorem
\ref{thmaspnw},
we shall denote for all
$x\in\dR$
\[
\wh{h}_{n}(x)=\frac{1}{n}\sum_{k=1}^nW_{k}(x)Y_{k}
\quad\mbox{and}\quad
\wh{g}_{n}(x)=\frac{1}{n}\sum_{k=1}^nW_{k}(x).
\]
As in~\cite{MR2448468}, we obtain from~(\ref{Sempar}) the decomposition
%
%e6.1 #&#
%e6.2 #&#
%
\begin{eqnarray}
\label{maindeco1}
n\wh{h}_{n}(x)&=&M_n(x)+P_n(x)+Q_n(x)+n\wh{g}_{n}(x)f(x), \\
\label{maindeco2}
n\wh{g}_{n}(x)&=&N_n(x)+R_n(x) +ng(\theta+x),
\end{eqnarray}
where
%
%e6.3 #&#
%e6.4 #&#
%
\begin{eqnarray}
\label{DefMn}
M_n(x) &=& \sum_{k=1}^{n}W_{k}(x)\veps_{k},
\\
\label{DefNn}
N_n(x) &=& \sum_{k=1}^{n}W_{k}(x)-\dE[W_{k}(x)|\cF_{k-1}]
\end{eqnarray}
and
%
%e6.5 #&#
%e6.6 #&#
%e6.7 #&#
%
\begin{eqnarray}
\label{DefPn}
P_{n}(x)&=&\sum_{k=1}^{n}W_{k}(x)\bigl(f(X_{k}-\wh{\theta}_{k-1})-f(x)\bigr),
\\
\label{DefQn}
Q_{n}(x)&=&\sum_{k=1}^{n}W_{k}(x)\bigl(f(X_{k}-\theta)-f(X_{k}-\wh{\theta
}_{k-1})\bigr),\\
\label{DefRn}
R_n(x) &=& \sum_{k=1}^{n}\bigl(\dE[W_{k}(x)|\cF_{k-1}]-g(\theta+x)\bigr).
\end{eqnarray}
On the one hand,
\[
\dE[W_{n}(x)|\cF_{n-1}]=\int_{\dR}\frac{1}{h_{n}}K\biggl(\frac
{x_{n}-\widehat{\theta}_{n-1}-x}{h_{n}}\biggr)g(x_n) \,dx_n.
\]
After the change of variables $z=h_n^{-1}(x_{n}-\widehat{\theta
}_{n-1}-x)$, as the density function~$g$ is
continuous, twice differentiable with bounded derivatives, we infer
from the Taylor formula
that
%
%e6.8 #&#
%
\begin{eqnarray}
\label{ExpW}
&&\dE[W_{n}(x)|\cF_{n-1}] \nonumber\\
&&\qquad= \int_{\dR}K(z) g(\wh{\theta
}_{n-1}+x+h_{n}z) \,dz\nonumber\\
&&\qquad=\int_{\dR}K(z)\biggl( g(\wh{\theta}_{n-1}+x)+h_{n}zg^{\prime
}(\wh{\theta}_{n-1}+x) \\
&&\hspace*{52pt}\qquad\quad{}+\frac{h_{n}^{2}z^{2}}{2}
g^{\prime\prime}(\wh{\theta}_{n-1}+x+h_n z\xi)\biggr)\,dz \nonumber\\
&&\qquad=g(\wh{\theta}_{n-1}+x)+\frac{h_{n}^{2}}{2}\int_{\dR
}z^{2}K(z)g^{\prime\prime}(\wh{\theta}_{n-1}+x+h_n z\xi)\,dz,\nonumber
\end{eqnarray}
where $0<\xi<1$. Consequently, for all $n\geq1$,
%
%e6.9 #&#
%
\begin{equation}
\label{sharpExpW1}
\bigl|\dE[W_{n}(x)|\cF_{n-1}] -g(\wh{\theta}_{n-1}+x)\bigr|\leq M_{g}\tau
^{2}h_{n}^{2} \qquad\mbox{a.s.},
\end{equation}
where
$M_{g}={\sup_{x \in\dR}} |g^{\prime\prime}(x)|$
and
\[
\tau^{2}=\frac{1}{2} \int_{\dR}x^{2}K(x)\,dx.
\]
The continuity of $g$ together with the fact that $\wh{\theta}_{n}$
converges to $\theta$ a.s. leads to
%
%e6.10 #&#
%
\begin{equation}
\label{Cvgg}
\lim_{n \rightarrow\infty}\frac{1}{n}\sum_{k=1}^{n} \dE
[W_{k}(x)|\cF_{k-1}] = g(\theta+x)
\qquad\mbox{a.s.},
\end{equation}
which immediately implies that for all $x \in\dR$
%
%e6.11 #&#
%
\begin{equation}
\label{CvgRn}
R_n(x)=o(n) \qquad\mbox{a.s.}\vadjust{\goodbreak}
\end{equation}
On the other hand, $(N_n(x))$ is a square integrable martingale
difference sequence with predictable quadratic variation
given by
\begin{eqnarray*}
\langle N(x) \rangle_n&=&\sum_{k=1}^n \dE\bigl[\bigl(N_k(x) - N_{k-1}(x)\bigr)^2|\cF_{k-1}\bigr] \\
&=&\sum_{k=1}^n \dE[W_{k}^2(x)|\cF_{k-1}]- \dE^2[W_{k}(x)|\cF_{k-1}].
\end{eqnarray*}
It follows from the same calculation as in~(\ref{ExpW}) that
\begin{eqnarray*}
\dE[W_{n}^2(x)|\cF_{n-1}] &=& \frac{1}{h_n} \int_{\dR}K^2(z) g(\wh
{\theta}_{n-1}+x+h_{n}z) \,dz \\
&=&\frac{\nu^2}{h_n} g(\wh{\theta}_{n-1}+x)+\frac{h_{n}}{2}\int
_{\dR}z^{2}K^2(z)g^{\prime\prime}(\wh{\theta}_{n-1}+x+h_n z\xi)\,dz,
\end{eqnarray*}
where $0<\xi<1$, which leads to
%
%e6.12 #&#
%
\begin{equation}
\label{sharpExpW2}
\biggl|\dE[W_{n}^2(x)|\cF_{n-1}] -\frac{\nu^2}{h_n}g(\wh{\theta
}_{n-1}+x)\biggr|\leq M_{g}\mu^{2}h_{n} \qquad\mbox{a.s.}
\end{equation}
with
\[
\nu^{2}= \int_{\dR}K^2(x)\,dx
\quad\mbox{and}\quad
\mu^{2}= \frac{1}{2}\int_{\dR}x^2K^2(x)\,dx.
\]
Hence, since
\[
\lim_{n\rightarrow\infty} \frac{1}{n^{1+\alpha}}\sum_{k=1}^n
h_{k}^{-1} = \frac{1}{1+\alpha}
\]
we deduce from~(\ref{sharpExpW1}) and~(\ref{sharpExpW2}) together
with the Toeplitz lemma
and the almost sure convergence of $g(\wh{\theta}_{n}+x)$
to $g(\theta+x)$ that
%
%e6.13 #&#
%
\begin{equation}
\label{CvgBracketNn}
\lim_{n\rightarrow\infty} \frac{\langle N(x) \rangle_n}{n^{1+\alpha}} =
\frac{\nu^2 g(\theta+x) }{1+\alpha}
\qquad\mbox{a.s.}
\end{equation}
Consequently, we obtain from the strong law of large numbers for
martingales given, for example, by Theorem 1.3.15 of~\cite{MR1485774}
that for any $\gamma> 0$,
$(N_n(x))^2=o(n^{1+\alpha}(\log n)^{1+\gamma}) $ a.s.
which ensures that, for all $x \in\dR$
%
%e6.14 #&#
%
\begin{equation}
\label{CvgNn}
N_n(x)=o(n) \qquad\mbox{a.s.}
\end{equation}
Therefore, it follows from~(\ref{maindeco2}),~(\ref{CvgRn}) and
(\ref{CvgNn}) that
for all $x \in\dR$
%
%e6.15 #&#
%
\begin{equation}
\label{Cvggchapn}
\lim_{n\rightarrow\infty}
\wh{g}_{n}(x)=g(\theta+x) \qquad\mbox{a.s.}
\end{equation}
Moreover, the kernel $K$ is compactly supported
which means that one can find a positive constant $A$ such
that $K$ vanishes outside the interval $[-A, A]$. Thus, for all $n \geq
1$ and all $x \in\dR$,
\[
W_n(x)=\frac{1}{h_{n}}K\biggl(\frac{X_{n}-\widehat{\theta
}_{n-1}-x}{h_{n}}\biggr)\rI_{\{|X_{n}-\widehat{\theta}_{n-1}-x| \leq
A h_n \}}.\vadjust{\goodbreak}
\]
In addition, the function~$f$ is Lipschitz, so there exists a positive
constant~$C_f$ such that for all $n \geq1$
\[
|f(X_{n}-\widehat{\theta}_{n-1})-f(x)| \leq C_f |X_{n}-\widehat
{\theta}_{n-1}-x|.
\]
Consequently, we obtain from~(\ref{DefPn}) that for all $x \in\dR$
%
%e6.16 #&#
%
\begin{eqnarray}\label{MajPn}
|P_n(x)| &\leq& C_f \sum_{k=1}^{n}W_{k}(x) |X_{k}-\widehat{\theta
}_{k-1}-x| \nonumber\\[-8pt]\\[-8pt]
&\leq& A C_f \sum_{k=1}^{n}h_k W_{k}(x).
\nonumber
\end{eqnarray}
Hence, it follows from convergence~(\ref{Cvgg}) together with~(\ref{CvgNn}) and~(\ref{MajPn}) that
for all $x \in\dR$
%
%e6.17 #&#
%
\begin{equation}
\label{CvgPn}
P_n(x)= o(n) \qquad\mbox{a.s.}
\end{equation}
Furthermore, we obtain from~(\ref{DefQn}) that for all $x \in\dR$
%
%e6.18 #&#
%
\begin{equation}\label{MQn}
|Q_n(x)| \leq C_f \sum_{k=1}^{n}W_{k}(x) |\widehat{\theta
}_{k-1}-\theta|.
\end{equation}
Then, it follows from the Cauchy--Schwarz inequality that
%
%e6.19 #&#
%
\begin{equation}\label{MajQn}
Q_n^2(x) \leq C_f^2 \sum_{k=1}^{n}W_{k}^2(x) \sum_{k=1}^{n}|\widehat
{\theta}_{k-1}-\theta|^2.
\end{equation}
We can split the first sum at the right-hand side of~(\ref{MajQn})
into two terms,
\[
\sum_{k=1}^{n}W_{k}^2(x)= I_n(x)+J_n(x),
\]
where
\begin{eqnarray*}
I_n(x)&=&\sum_{k=1}^{n}W_{k}^2(x)- \dE[W_{k}^2(x)|\cF_{k-1}],\\
J_n(x)&=&\sum_{k=1}^{n} \dE[W_{k}^2(x)|\cF_{k-1}].
\end{eqnarray*}
Following the same lines as in the proof of~(\ref{CvgNn}), it is not
hard to see that
\[
I_n(x)=o(n^{1+\alpha}) \qquad\mbox{a.s.}
\]
We also deduce from convergence~(\ref{CvgBracketNn}) that
\[
J_n(x)=O(n^{1+\alpha}) \qquad\mbox{a.s.}
\]
Consequently, we obtain that for all $x \in\dR$
%
%e6.20 #&#
%
\begin{equation}
\label{MajsumW}
\sum_{k=1}^{n}W_{k}^2(x)= O(n^{1+\alpha}) \qquad\mbox{a.s.}\vadjust{\goodbreak}
\end{equation}
Therefore, we infer from the quadratic strong law given by~(\ref{qslrm}) together
with~(\ref{MajQn}) and~(\ref{MajsumW}) that $Q_n^2(x)= O(n^{1+\alpha
} \log n)$ a.s. which
implies that for all $x \in\dR$
%
%e6.21 #&#
%
\begin{equation}
\label{CvgQn}
Q_n(x)= o(n) \qquad\mbox{a.s.}
\end{equation}
It now remains to study the asymptotic behavior of $M_n(x)$ given by
(\ref{DefMn}).
As~$(X_n)$ and $(\veps_n)$ are two independent sequences of
independent and identically distributed random variables,
$(M_n(x))$ is a square integrable martingale difference sequence with
predictable quadratic variation
given by
\begin{eqnarray*}
\langle M(x) \rangle_n&=&\sum_{k=1}^n \dE\bigl[\bigl(M_k(x) - M_{k-1}(x)\bigr)^2|\cF_{k-1}\bigr] \\
&=&\sigma^2 \sum_{k=1}^n \dE[W_{k}^2(x)|\cF_{k-1}].
\end{eqnarray*}
Then, it follows from convergence~(\ref{CvgBracketNn}) that
%
%e6.22 #&#
%
\begin{equation}
\label{CvgBracketMn}
\lim_{n\rightarrow\infty} \frac{\langle M(x)\rangle_n}{n^{1+\alpha}} =
\frac{\sigma^2\nu^2 g(\theta+x) }{1+\alpha}
\qquad\mbox{a.s.}
\end{equation}
Consequently, we obtain from the strong law of large numbers for
martingales that
for any $\gamma> 0$,
$(M_n(x))^2=o(n^{1+\alpha}(\log n)^{1+\gamma}) $ a.s.
which leads to
%
%e6.23 #&#
%
\begin{equation}
\label{CvgMn}
M_n(x)=o(n) \qquad\mbox{a.s.}
\end{equation}
Therefore, we deduce from~(\ref{maindeco1}) and~(\ref{Cvggchapn})
together with
the conjunction of~(\ref{CvgPn}),~(\ref{CvgQn}) and~(\ref{CvgMn})
that for all $x \in\dR$
%
%e6.24 #&#
%
\begin{equation}
\label{Cvghchapn}
\lim_{n\rightarrow\infty}
\wh{h}_{n}(x)=f(x)g(\theta+x) \qquad\mbox{a.s.}
\end{equation}
Finally, we can conclude from the identity
%
%e6.25 #&#
%
\begin{equation}
\label{Idfn}
\wh{f}_{n}(x)= \frac{\wh{h}_{n}(x)+\wh{h}_{n}(-x)}{\wh
{g}_{n}(x)+\wh{g}_{n}(-x)}
\end{equation}
and the parity of the function~$f$ that, for all $x \in\dR$ such that
$|x|\leq{1/2}$,
%
%e6.26 #&#
%
\begin{equation}
\lim_{n\rightarrow\infty}
\wh{f}_{n}(x)=f(x) \qquad\mbox{a.s.}
\end{equation}

%s6.2 #&#
\subsection{\texorpdfstring{Proof of Theorem \protect\ref{thmcltnw}}{Proof of Theorem 3.2}}\label{sec62}

We shall now proceed to the proof of the asymptotic normality of $\wh{f}_{n}$.
It follows from~(\ref{maindeco1}),~(\ref{maindeco2}) and~(\ref{Idfn}) that for all $x \in\dR$
%
%e6.27 #&#
%
\begin{equation}
\label{deltafn}
\wh{f}_{n}(x)-f(x)=\frac{\cM_n(x)+\cP_n(x)+\cQ_n(x)}{n\cG_n(x)},
\end{equation}
where $\cG_n(x)=\wh{g}_{n}(x)+\wh{g}_{n}(-x)$ and
\begin{eqnarray*}
\cM_n(x)&=&M_n(x)+M_n(-x), \\
\cP_n(x)&=&P_n(x)+P_n(-x), \\
\cQ_n(x)&=&Q_n(x)+Q_n(-x)
\end{eqnarray*}
with $M_n(x)$, $P_n(x)$ and $Q_n(x)$ given by
(\ref{DefMn}),~(\ref{DefPn}) and~(\ref{DefQn}), respectively. We
already saw from~(\ref{Cvggchapn}) that for all $x \in\dR$
%
%e6.28 #&#
%
\begin{equation}
\label{CvgGn}
\lim_{n\rightarrow\infty}
\cG_n(x)=g(\theta+x)+g(\theta-x) \qquad\mbox{a.s.}
\end{equation}
In order to establish the asymptotic normality, it is now necessary to
be more precise in the almost sure rates
of convergence given in~(\ref{CvgPn}) and~(\ref{CvgQn}). It follows
from~(\ref{MajPn}) that for all $x \in\dR$
%
%e6.29 #&#
%
\begin{equation}
\label{MajPPn}
|P_n(x)| \leq A C_f \bigl(L_n(x)+ \Lambda_n(x)\bigr),
\end{equation}
where
\begin{eqnarray*}
L_n(x)&=&\sum_{k=1}^{n}h_k \bigl(W_{k}(x) - \dE[W_{k}(x) | \cF_{k-1}]\bigr),\\
\Lambda_n(x)&=&\sum_{k=1}^{n}h_k\dE[W_{k}(x)| \cF_{k-1}].
\end{eqnarray*}
On the one hand, we infer from~(\ref{sharpExpW1}) that
%
%e6.30 #&#
%
\begin{equation}
\label{MajLambdan}
\Lambda_n(x)=
O\Biggl(\sum_{k=1}^{n}h_k\Biggr)
=O(n^{1-\alpha})
\qquad\mbox{a.s.}
\end{equation}
On the other hand, $(L_n(x))$ is a square integrable martingale
difference sequence with predictable quadratic variation
given by
\[
\langle L(x)\rangle_n= \sum_{k=1}^n h_k^2\bigl(\dE[W_{k}^2(x)|\cF_{k-1}] - \dE
^2[W_{k}(x)| \cF_{k-1}]\bigr).
\]
We deduce from~(\ref{sharpExpW1}) and~(\ref{sharpExpW2}) together
with the Toeplitz lemma that
%
%e6.31 #&#
%
\begin{equation}
\label{CvgBracketLn}
\lim_{n\rightarrow\infty} \frac{\langle L(x)\rangle_n}{n^{1-\alpha}} =
\frac{\nu^2 g(\theta+x) }{1-\alpha}
\qquad\mbox{a.s.}
\end{equation}
Consequently, we obtain from the strong law of large numbers for
martingales that for any $\gamma> 0$,
$(L_n(x))^2=o(n^{1-\alpha}(\log n)^{1+\gamma}) $ a.s.
which clearly implies that $(L_n(x))^2=o(n^{1+\alpha}) $ a.s.
Therefore, we find from~(\ref{MajPPn}) and~(\ref{MajLambdan}) that,
as soon as $\alpha>1/3$,
\[
(P_n(x))^2 =O(n^{2-2\alpha})+o(n^{1+\alpha})= o(n^{1+\alpha})
\qquad\mbox{a.s.},
\]
which immediately leads to
%
%e6.32 #&#
%
\begin{equation}
\label{RescalPn}
(\cP_n(x))^2 = o(n^{1+\alpha}) \qquad\mbox{a.s.}\vadjust{\goodbreak}
\end{equation}
Proceeding as in the proof of~(\ref{RescalPn}), we obtain from~(\ref{MQn}) that for all $x \in\dR$
%
%e6.33 #&#
%
\begin{equation}
\label{MajQQn}
|Q_n(x)| \leq C_f \bigl(S_n(x)+ \Sigma_n(x)\bigr),
\end{equation}
where
\begin{eqnarray*}
S_n(x)&=&\sum_{k=1}^{n}\ell_k \bigl(W_{k}(x) - \dE[W_{k}(x) | \cF
_{k-1}]\bigr),\\
\Sigma_n(x)&=&\sum_{k=1}^{n}\ell_k\dE[W_{k}(x)| \cF_{k-1}]
\end{eqnarray*}
with $\ell_n= |\widehat{\theta}_{n-1}-\theta|$. We deduce from
(\ref{sharpExpW1})
together with the Cauchy--Schwarz inequality and the quadratic strong
law given by~(\ref{qslrm})
that
%
%e6.34 #&#
%
\begin{equation}
\label{MajSigman}
\Sigma_n(x)=
O\Biggl(\sum_{k=1}^{n}\ell_k\Biggr)
=O\bigl(\sqrt{n \log n}\bigr)
\qquad\mbox{a.s.}
\end{equation}
In addition, it follows from~(\ref{sharpExpW2}) that
$(S_n(x))$ is a square integrable martingale difference sequence with
predictable quadratic variation satisfying
\[
\langle S(x)\rangle_n=O(n^{\alpha} \log n) \qquad\mbox{a.s.}
\]
Consequently, we obtain from the strong law of large numbers for
martingales that for any $\gamma> 0$,
$(S_n(x))^2=o(n^{\alpha}(\log n)^{2+\gamma}) $ a.s. so
$(S_n(x))^2=o(n^{1+\alpha}) $ a.s. Hence, we find from~(\ref{MajQQn})
and~(\ref{MajSigman}) that
\[
(Q_n(x))^2 =O(n \log n)+o(n^{1+\alpha})= o(n^{1+\alpha})
\qquad\mbox{a.s.},
\]
which obviously implies
%
%e6.35 #&#
%
\begin{equation}
\label{RescalQn}
(\cQ_n(x))^2 = o(n^{1+\alpha}) \qquad\mbox{a.s.}
\end{equation}
It remains to establish the asymptotic behavior of the dominating term~$\cM_n(x)$.
We already saw that $(M_n(x))$ is a square integrable martingale
difference sequence.
Consequently, $(\cM_n(x))$ is also a square integrable martingale
difference sequence
with predictable quadratic variation given by
\[
\langle \cM(x)\rangle_n=\sigma^2 \sum_{k=1}^n \dE\bigl[\bigl(W_{k}(x)+W_{k}(-x)\bigr)^2|\cF_{k-1}\bigr].
\]
Hence, it is necessary to evaluate the cross-term $\dE
[W_{n}(x)W_{n}(-x)|\cF_{n-1}]$.
It follows from the same calculation as in~(\ref{ExpW}) that
\begin{eqnarray*}
&&
\dE[W_{n}(x)W_{n}(-x)|\cF_{n-1}] \\[-2pt]
&&\qquad= \frac{1}{h_n} \int_{\dR
}K(z)K(z+2h_n^{-1}x) g(\wh{\theta}_{n-1}+x+h_{n}z) \,dz \\[-2pt]
&&\qquad=\frac{1}{h_n} g(\wh{\theta}_{n-1}+x)I_n(x)+g^{\prime}(\wh
{\theta}_{n-1}+x)J_n(x)\\[-2pt]
&&\qquad\quad{}+
\frac{h_{n}}{2}\int_{\dR}z^{2}K(z)K(z+2h_n^{-1}x)g^{\prime\prime
}(\wh{\theta}_{n-1}+x+h_n z\xi)\,dz
\end{eqnarray*}
with $0<\xi<1$. Consequently, we obtain that
\begin{eqnarray*}
&&\biggl|\dE[W_{n}(x)W_{n}(-x)|\cF_{n-1}] - \frac{1}{h_n} g(\wh
{\theta}_{n-1} + x)I_n(x) - g^{\prime}(\wh{\theta}_{n-1} + x)J_n(x)
\biggr|\\
&&\qquad \leq M_{g}H_n(x)h_{n} \qquad\mbox{a.s.},
\end{eqnarray*}
where
\begin{eqnarray*}
I_n(x)&=&\int_{\dR}K(z)K(z+2h_n^{-1}x)\,dz, \\[-2pt]
J_n(x)&=&\int_{\dR}zK(z)K(z+2h_n^{-1}x)\,dz, \\[-2pt]
H_n(x)&=&\int_{\dR}z^2 K(z)K(z+2h_n^{-1}x)\,dz.
\end{eqnarray*}
However, as the kernel $K$ is compactly supported, we have for all
$x\in\dR$ with $x \neq0$,
\[
\lim_{n \rightarrow\infty} K(z+2h_n^{-1}x)=0.
\]
Then, we deduce from the Lebesgue dominated convergence theorem
that all the three integrals $I_n(x)$, $J_n(x)$ and $H_n(x)$ tend to zero
as $n$ goes to infinity, which implies that for all $x\in\dR$ with $x
\neq0$,
%
%e6.36 #&#
%
\begin{equation}
\label{sharpExpWcross}\qquad
\sum_{k=1}^n \dE[W_{k}(x)W_{k}(-x)|\cF_{k-1}]= o\Biggl( \sum
_{k=1}^n h_{k}^{-1}\Biggr)=o(n^{1+\alpha}) \qquad\mbox{a.s.}
\end{equation}
Therefore, we find from~(\ref{CvgBracketMn}) together with~(\ref{sharpExpWcross}) that for all $x\in\dR$ with $x \neq0$,
%
%e6.37 #&#
%
\begin{equation}
\label{CvgBracketcalMn}
\lim_{n\rightarrow\infty} \frac{\langle \cM(x)\rangle_n}{n^{1+\alpha}} =
\frac{\sigma^2\nu^2}{1+\alpha}\bigl(g(\theta+x)+g(\theta-x)\bigr)
\qquad\mbox{a.s.}
\end{equation}
If $x=0$, it immediately follows from~(\ref{CvgBracketMn})
%
%e6.38 #&#
%
\begin{equation}
\label{CvgBracketcalMnzero}
\lim_{n\rightarrow\infty} \frac{\langle \cM(0)\rangle_n}{n^{1+\alpha}} =
\frac{4\sigma^2\nu^2 g(\theta)}{1+\alpha}
\qquad\mbox{a.s.}
\end{equation}
Furthermore, it is not hard to see that the Lindeberg condition is
satisfied. As a matter of fact,
we have assumed that the sequence $(\varepsilon_n)$ has a finite
moment of order $a>2$. If we denote
$\Delta\cM_n(x)=\cM_n(x) - \cM_{n-1}(x)$, we have
\[
\dE[| \Delta\cM_n(x) |^a|\cF_{n-1}] =\dE[|\varepsilon_n|^a] \dE
[| W_n(x)- W_n(-x) |^a|\cF_{n-1}],
\]
which implies that
\[
\dE[| \Delta\cM_n(x) |^a|\cF_{n-1}]
\leq2^{a-1} \dE[|\varepsilon_n|^a] \dE[ W_n^a(x) + W_n^a(-x) |\cF_{n-1}].
\]
However, it follows from the same calculation as in~(\ref{ExpW}) that
%
%e6.39 #&#
%
\begin{equation}
\label{sharpExpWa}\qquad
\sum_{k=1}^n \dE[W_{k}^a(x)|\cF_{k-1}]= O\Biggl( \sum_{k=1}^n
h_{k}^{1-a}\Biggr)=O\bigl(n^{1+\alpha(a-1)}\bigr) \qquad\mbox{a.s.}\vadjust{\goodbreak}
\end{equation}
In addition, for all $\varepsilon>0$,
\begin{eqnarray*}
&&\frac{1}{n^{1+\alpha}}\sum_{k=1}^{n}\dE\bigl[(\Delta\cM_k(x))^{2}\rI
_{|\Delta\cM_k(x)|\geq{\veps\sqrt{n^{1+\alpha}}}}|\cF_{k-1}\bigr]\\
&&\qquad\leq
\frac{1}{\veps^{a-2} n^b}\sum_{k=1}^{n}\dE[| \Delta\cM_k(x)
|^a|\cF_{k-1}],
\end{eqnarray*}
where $b=a(1+\alpha)/2$. Consequently, it follows from~(\ref{sharpExpWa})
that for all $\varepsilon>0$,
\[
\frac{1}{n^{1+\alpha}}\sum_{k=1}^{n}\dE\bigl[(\Delta\cM_k(x))^{2}\rI
_{|\Delta\cM_k(x)|\geq{\veps\sqrt{n^{1+\alpha}}}}|\cF_{k-1}\bigr]=O(n^c)
\qquad\mbox{a.s.},
\]
where $c=(2-a)(1-\alpha)/2$. As $c<0$, the Lindeberg condition is
clearly satisfied. We can conclude
from the central limit theorem for martingales given, for example, by
Corollary 2.1.10 of~\cite{MR1485774}
that for all
$x \in\dR$ with $x \neq0$,
%
%e6.40 #&#
%
\begin{equation}
\label{cltMn}
\frac{\cM_{n}(x)}{\sqrt{n^{1+\alpha}}}\liml\cN\biggl(0,\frac
{\sigma^{2}\nu^2}{1+\alpha}\bigl(g(\theta+x)+g(\theta-x)\bigr)\biggr),
\end{equation}
while, for $x=0$,
%
%e6.41 #&#
%
\begin{equation}
\label{cltMnzero}
\frac{\cM_{n}(0)}{\sqrt{n^{1+\alpha}}}\liml\cN\biggl(0,\frac
{4\sigma^{2}\nu^2}{1+\alpha}g(\theta)\biggr).
\end{equation}
Finally, it follows from~(\ref{deltafn}) and~(\ref{CvgGn}) together
with~(\ref{RescalPn}),~(\ref{RescalQn}),~(\ref{cltMn}), (\ref{cltMnzero})
and the Slutsky lemma that, for all $x \in\dR$ such that $|x|\leq
{1/2}$ with $x \neq0$,
\[
\sqrt{nh_n}\bigl(\wh{f}_{n}(x)-f(x)\bigr) \liml\cN\biggl(0,\frac{\sigma
^{2}\nu^2}{(1+\alpha)(g(\theta+x)+g(\theta-x))}\biggr),
\]
while, for $x=0$,
\[
\sqrt{nh_n}\bigl(\wh{f}_{n}(0)-f(0)\bigr) \liml\cN\biggl(0,\frac{\sigma
^{2}\nu^2}{(1+\alpha)g(\theta)}\biggr),
\]
which completes the proof of Theorem~\ref{thmcltnw}.

\section*{Acknowledgments}

The authors would like to thank the Associate Editor and the two
anonymous reviewers for their suggestions and constructive comments
which helped to improve the paper substantially.

%suskaldyti doi

% imsref loaded by lrinkeviciute, 2012-03-14 08:55:14
%

\printaddresses

\end{document}